\documentclass[twocolumn]{autart} 

\usepackage[]{graphicx}
\usepackage[usenames,dvipsnames]{color}
\usepackage{epstopdf}
\usepackage{booktabs}

\usepackage{numcompress}

\usepackage{amsmath,amssymb,mathrsfs,dsfont,mathdots}
\usepackage{algorithmic}
\usepackage[section]{algorithm}

\newtheorem{definition}[thm]{Definition}
\newtheorem{corollary}[thm]{Corollary}
\newtheorem{proposition}[thm]{Proposition}
\newtheorem{problem}[thm]{Problem}
\newtheorem{remark}[thm]{Remark}
\newtheorem{assumption}[thm]{Assumption}

\newcommand*\mc[0]{\mathcal}                                                  		

\newcommand{\until}[1]{\{1,\dots,#1\}}

\newcommand*\diag[0]{\mbox{diag}} 

\usepackage[prependcaption,colorinlistoftodos]{todonotes}

\newcommand\oprocendsymbol{\hbox{$\square$}}
\newcommand\oprocend{\relax\ifmmode\else\unskip\hfill\fi\oprocendsymbol}

\newcommand{\real}[0]{\mathbb R}

\DeclareSymbolFont{bbold}{U}{bbold}{m}{n}
\DeclareSymbolFontAlphabet{\mathbbold}{bbold}

\newcommand{\vect}[1]{\mathbbold{#1}}
\newcommand{\vectorones}[1][]{\vect{1}_{#1}}
\newcommand{\vectorzeros}[1][]{\vect{0}_{#1}}
\newcommand{\col}{\text{col}}

\def\begequ{\begin{equation}}
\def\endequ{\end{equation}}

\begin{document}

	\begin{frontmatter}
\title{Robustness of Distributed Averaging Control in Power Systems: Time Delays \& Dynamic Communication Topology}

 \thanks[footnoteinfo]{The project leading to this manuscript has received funding from the European Union's Horizon 2020 research and innovation programme under the Marie Sk\l odowska-Curie grant agreement No. 734832, the Ministry of Education and Science of Russian Federation (project 14.Z50.31.0031), ETH Z\"urich funds, the SNF Assistant Professor Energy Grant \#160573 and the Israel Science Foundation (Grant No. 1128/14).}
\author[Paestum]{Johannes Schiffer}\ead{j.schiffer@leeds.ac.uk}    
\quad \quad
\author[Rome]{Florian D\"orfler}\ead{dorfler@control.ee.ethz.ch}    
\quad \quad
\author[fri,itmo]{Emilia Fridman}\ead{emilia@eng.tau.ac.il}  

\address[Paestum]{School of Electronic \& Electrical Engineering, University of Leeds, United Kingdom, LS2 9JT}             
\address[Rome]{Automatic Control Laboratory, Swiss Federal Institute of Technology (ETH) Z\"urich, Switzerland, 8092}  
\address[fri]{Tel Aviv University, Tel Aviv 69978, Israel}
\address[itmo]{ITMO University, 49 Avenue Kronverkskiy, 197101 Saint Petersburg, Russia}

 
 \begin{abstract}
Distributed averaging-based integral (DAI) controllers are becoming increasingly popular in power system applications. The literature has thus far primarily focused on disturbance rejection, steady-state optimality and adaption to complex physical system models without considering uncertainties on the cyber and communication layer nor their effect on robustness and performance.
In this paper, we derive sufficient delay-dependent conditions for robust stability of a secondary-frequency-DAI-controlled power system with respect to heterogeneous communication delays, link failures and packet losses. Our analysis takes into account both constant as well as fast-varying delays, and it is based on a common strictly decreasing Lyapunov-Krasovskii functional.
The conditions illustrate an inherent trade-off between robustness and performance of DAI controllers. 
The effectiveness and tightness of our stability certificates are illustrated via a numerical example based on Kundur's four-machine-two-area test system.
 \end{abstract}



\end{frontmatter}



\section{Introduction}
\label{Section: Introduction}

\subsection{Motivation}

Power systems worldwide are currently experiencing drastic changes and challenges. One of the main driving factors for this development is the increasing penetration of distributed and volatile renewable generation interfaced to the network with power electronics accompanied by a reduction in synchronous generation. This results in power systems being operated under more and more stressed conditions \cite{WW-KE-GB-JK:15}. In order to successfully cope with these changes, the control and operation paradigms of today's power systems have to be adjusted to the new conditions. The increasing complexity in terms of network dynamics and number of active network elements renders centralized and inflexible approaches inappropriate creating a clear need for robust and distributed solutions with plug-and-play capabilities \cite{GS-NH-JPL-CM-AD-DP:15}. The latter approaches require a combination of advanced control techniques with adequate communication technologies. 

Multi-agent systems (MAS) represent a promising framework to address these challenges \cite{SM-ED-VC-AD-NH-FP-FT:07}. A popular distributed control strategy for MAS are distributed averaging-based integral (DAI) algorithms, also known as consensus filters \cite{freeman2006stability,ROS-AF-RM:07} that rely on  averaging of integral actions through a communication network.
The distributed character of this type of protocol has the advantage that no central computation unit is needed and the individual agents, i.e., generation units, only have to exchange information with their neighbors \cite{AB-FL-AD:14}. 

One of the most relevant control applications in power systems is frequency control which is typically divided into three hierarchical layers: primary, secondary and tertiary control \cite{PK:94}. In the present paper, we focus on secondary control which is tasked with the regulation of the frequency to a nominal value in an economically efficient way and subject to maintaining the net area power balance. 
The literature on secondary DAI frequency controllers addressing these tasks is reviewed in the following.%

\subsection{Literature Review on DAI Frequency Control}

DAI algorithms have been proposed previously to address the objectives of secondary frequency control in bulk power systems \cite{ST-MB-CDP:16,MA-HS-VDD-KHJ:13,JS-FD:15,NM-CDP:15} and also in microgrids (i.e., small-footprint power systems on the low and medium voltage level) \cite{JWSP-FD-FB:12u,AB-FL-AD:14,FD-JWSP-FB:14a}. They have been extended to achieve asymptotically optimal injections \cite{TS-CDP-AVDS:15,CZ-EM-FD:15}, and have also been adapted to increasingly complex physical system models \cite{CDP-NM-JS-FD:16,TS-CDP-AVDS:16}. The closed-loop DAI-controlled power system is a cyber-physical system whose stability and performance crucially relies on nearest-neighbor communication. Despite all recent advances, communication-based controllers (in power systems) are subject to considerable uncertainties such as message delays, message losses, and link failures \cite{QY-JAB-TCG:11,GS-NH-JPL-CM-AD-DP:15} that can severely reduce the performance -- or even affect the stability -- of the overall cyber-physical system. 
Such cyber-physical phenomena and uncertainties have not been considered thus far in DAI-controlled power systems.

For microgrids, the effect of communication delays on secondary controllers has been considered in \cite{SL-XW-PXL:15} for the case of a centralized PI controller, in \cite{CA-RC-DS-JMG:15} for a centralized PI controller with a Smith predictor as well as a  model predictive controller and in \cite{EAC-DW-JMG-JCV-TD-CS-PP:16} for a DAI-controlled microgrid with fixed communication topology. In all three papers, a small-signal (i.e., linearization-based) analysis of a model with constant delays is performed. 

In \cite{JL-HZ-XL-LZ:16,JL-HZ-XL-XY-WH:16} distributed control schemes for microgrids are proposed, and conditions for stability under time-varying delays as well as a dynamic communication topology are derived. However, both approaches are based on the pinning-based controllers requiring a master-slave architecture. Compared to the DAI controller in the present paper, this introduces an additional uncertainty as the leader may fail (see also Remark 1 in \cite{JL-HZ-XL-XY-WH:16}). In addition, the analysis in \cite{JL-HZ-XL-LZ:16,JL-HZ-XL-XY-WH:16} is restricted to the distributed control scheme on the cyber layer and neglects the physical dynamics. 
Moreover, the control in \cite{JL-HZ-XL-LZ:16} is limited to power sharing strategies, and secondary frequency regulation is not considered.

The delay robustness of alternative distributed secondary control strategies (based on primal-dual decomposition approaches) has been investigated for constant delays and a linearized power system model in \cite{XZ-AP:13,XZ-RK-MM-AP:16}.

\subsection{Contributions}
The present paper addresses both the cyber and the physical aspects of DAI frequency control by deriving conditions for robust stability of nonlinear DAI-controlled power systems under communication uncertainties. 
With regards to delays, we consider constant as well as fast-varying delays. The latter are a common phenomenon in sampled data networked control systems, due to digital control \cite{KL-EF:12,EF:14,EF:14b} and as the network access and transmission delays depend on the actual network conditions, e.g., in terms of congestion and channel quality \cite{JPH-PN-YX:07}. 
In addition to delays, in practical applications the topology of the communication network can be time-varying due to message losses and link failures \cite{ROS-RM:04,ROS-AF-RM:07,PL-YJ:08}. This can be modeled by a switching communication network \cite{ROS-RM:04,ROS-AF-RM:07}.  
Thus, the explicit consideration of communication uncertainties leads to a {\em switched nonlinear power system model with (time-varying) delays} the stability of which is investigated in this paper. 

More precisely, our main contributions are as follows.
First, we derive a strict Lyapunov function for a nominal DAI-controlled power system model without communication uncertainties, which may also be of independent interest.
Second, we extend this strict Lyapunov function to a common Lyapunov-Krasovskii functional (LKF) to provide sufficient delay-dependent conditions for robust stability of a DAI-controlled power system with dynamic communication topology as well as heterogeneous constant and fast-varying delays. Our stability conditions can be verified without exact knowledge of the operating state and reflect a fundamental trade-off between robustness and performance of DAI control.
Third and finally, we illustrate the effectiveness of the derived approach on a numerical benchmark example, namely Kundur's four-machine-two-area test system \cite[Example 12.6]{PK:94}.

The remainder of the paper is structured as follows. In Section~\ref{sec:prel} we recall some preliminaries on algebraic graph theory and introduce the power system model employed for the analysis. The DAI control is motivated and introduced in Section~\ref{sec: closed-loop model}, where we also derive a suitable error system. A strict Lyapunov function for the closed-loop DAI-controlled power system is derived in Section~\ref{sec:stabNom}. Based on this Lyapunov function, we then construct a common LKF for DAI-controlled power systems with constant and fast-varying delays in Section~\ref{sec:stabDel}.
A numerical example is provided in Section~\ref{sec: sim}. The paper is concluded with a brief summary and outlook on future work in Section~\ref{sec:Conclusions}.

\textbf{Notation.}
We define the sets $\real_{\geq0}:=\{x \in \real|x\geq0 \}$, $\real_{>0}:=\{x \in \real|x>0 \}$ and $\real_{<0}:=\{x \in \real|x<0 \}.$
For a set $\mathcal V,$ $|\mathcal V|$ denotes its cardinality and $[\mathcal V]^k$ denotes the set of all subsets of $\mathcal V$ that contain $k$ elements.
Let $x=\col(x_i)\in\real^n$ denote a vector with entries $x_i$ for $i \in \until n$, 
$\vectorzeros[n]$ the zero vector, $\vectorones[n]$ the vector with all entries equal to one, $I_{n}$ the $n\times n$ identity matrix, $\vectorzeros[n\times n]$ the $n\times n$ matrix with all entries equal to zero and $\diag(a_i)$ an
$n\times n$ diagonal matrix with diagonal entries $a_i\in\real.$ Likewise, $A=\text{blkdiag}(A_i)$ denotes a block-diagonal matrix with block-diagonal matrix entries $A_i.$
For $A\in\real^{n\times n},$ $A>0$ means that $A$ is symmetric positive definite. 
The elements below the diagonal of a symmetric matrix are denoted by $\ast.$
We denote by \mbox{$W[-h,0],$ $h\in\real_{>0},$} the Banach space of absolutely continuous functions $\phi:[-h,0]\to\real^n,$ $h\in\real_{>0},$ with $\dot \phi\in L_2(-h,0)^n$ and with the norm $\|\phi \|_W=\max_{\theta\in[a,b]}|\phi(\theta)|+\left({\int_{-h}^0{\dot\phi }^2d\theta}\right)^{0.5}$.
Also, $\nabla f$ denotes the gradient of a function $f:\real^n \to \real.$

\section{Preliminaries}
\label{sec:prel}
\subsection{Algebraic Graph Theory}
An undirected graph of order $n$ is a tuple \mbox{$\mathcal G =(\mathcal N, \mathcal E),$} where $\mathcal{N} = \{1, \ldots , n\}$ is the set of nodes and $\mathcal E \subseteq [\mathcal N]^2,$ 
$\mathcal E=\{e_1,\ldots,e_m\},$ is the set of undirected edges, i.e., the elements of $\mathcal E$ are subsets of $\mathcal N$ that contain two elements.
In the context of the present work, each node in the graph represents a generation unit.  
The adjacency matrix $\mathcal A\in \real^{|\mathcal N|\times |\mathcal N|}$ has entries $a_{ik}=a_{ki}=1$ if an edge between $i$ and $k$ exists and $a_{ik}=0$ otherwise. The degree of a node $i$ is defined as $d_i=\sum_{k=1}^{|\mathcal N|} a_{ik}.$ 
The Laplacian matrix of an undirected graph is given by 
$\mc L=\mathcal D - \mathcal A,$
where $\mathcal D=\diag(d_i)\in\real^{|\mathcal N|\times |\mathcal N|}.$ 
An ordered sequence of nodes such that any pair of consecutive nodes in the sequence is connected by an edge is called a path.
A graph $\mathcal G$ is called connected if for all pairs $\{i,k\} \in  [\mathcal N]^2$ there exists a path from $i$ to $k$. 
The Laplacian matrix $\mathcal L$ of an undirected graph is positive semidefinite with a simple zero eigenvalue if and only if the graph is connected. 
The corresponding right eigenvector to this simple zero eigenvalue is $\vectorones[n]$, i.e., $\mathcal L\vectorones[n]=\vectorzeros[n]$ \cite{CDG-GFR:01}. We refer the reader to \cite{RD:00,CDG-GFR:01} for further information on graph theory.

\subsection{Power Network Model}
We consider a Kron-reduced \cite{FD-FB:11d,PK:94} power system model composed of $n\geq1$ nodes. The set of network nodes is denoted by $\mc N=\{1,\ldots,n\}.$
Following standard practice, we make the following assumptions: the line admittances are purely inductive and the voltage amplitudes $V_i\in\real_{>0}$ at all nodes $i\in\mc N$ are constant \cite{PK:94}.
To each node $i\in\mc N,$ we associate a phase angle $\theta_i:\real_{\geq0}\to\real$ and denote its time derivative by $\omega_i=\dot \theta_i,$ which represents the electrical frequency at the $i$-th node. Under the made assumptions, two nodes $i$ and $k$ are connected via a nonzero susceptance $B_{ik}\in\real_{<0}.$ 
If there is no line between $i$ and $k,$ then $B_{ik}=0.$
We denote by $\mc N_i=\{k\in\mc N\,|\, B_{ik}\neq0\}$ the set of neighboring nodes of the $i$-th node.
Furthermore, we assume that for all $\{i,k\}\in[\mc N]^2$ there exists an ordered sequence of nodes from $i$ to $k$ such that any pair of consecutive nodes in the sequence is connected by a power line represented by an admittance, i.e., the electrical network is connected.

We consider a heterogeneous generation pool consisting of rotational synchronous generators (SGs) and inverter-interfaced units. The former are the standard equipment in conventional power networks and mostly used to connect fossil-fueled generation to the network. Compared to this, most renewable and storage units are connected via inverters, i.e., power electronics equipment, to the grid. We assume that all inverter-interfaced units are fitted with droop control and power measurement filters, see \cite{QCZ-TH:13,JS:15-thesis}. This implies that their dynamics admit a mathematically equivalent representation to SGs \cite{JS:15-thesis,JS-RO-AA-JR-TS:14}. Hence, the dynamics of the generation unit at the $i$-th node, $i\in \mc N,$ considered in this paper is given by
\begin{equation}
	\begin{split}
\dot \theta_{i}& = \omega_{i},
	\\
	M_{i} \dot \omega_{i} & = - D_{i} (\omega_{i}-\omega^d) + P_{i}^d - G_{ii}V_i^2 +u_i - P_i,
\end{split}
\label{eq:model}%
\end{equation}
where $D_i\in\real_{>0}$ is the damping or (inverse) droop coefficient, $\omega^d\in\real_{>0}$ is the nominal frequency, $P_i^d\in\real$ is the active power setpoint and $G_{ii}V_i^2,$ $G_{ii}\in\real_{\geq0},$ represents the (constant active power) load at the $i$-th node\footnote{For constant voltage amplitudes, any constant power load can equivalently be represented by a constant impedance load, i.e., to any constant $P\in\real_{>0}$ and constant $V\in\real_{>0},$ there exists a constant $G\in\real_{>0},$ such that $P=GV^2.$ On larger time scales, the loads may not be constant but follow regular (e.g., daily) fluctuations that can be accurately represented by internal models in the controllers \cite{ST-MB-CDP:16,NM-CDP:15,RP-CS-RW:16}. For the time-scales of interest to us, a constant load model suffices and, accordingly, our controllers contain integrators.}. Furthermore, $u_i:\real_{\geq0}\to\real$ is a control input and $M_i\in\real_{>0}$ is the (virtual) inertia coefficient, which in case of an inverter-interfaced unit is given by
$
M_i=\tau_{P_i}D_i,
$
where $\tau_{P_i}\in\real_{>0}$ is the low-pass filter time constant of the power measurement filter, see \cite{JS-DG-JR-TS:13,JS:15-thesis}. Following standard practice in power systems, all parameters are assumed to be given in per unit \cite{PK:94}.
The active power flow $P_i:\real^n\to\real$ is given by
\[
P_i=\sum\nolimits_{k \in \mathcal N_i} |B_{ik}| V_{i} V_{k} \sin(\theta_{ik}),
\]
where we have introduced the short-hand $\theta_{ik}=\theta_{i}-\theta_{k}.$
For a detailed modeling of the system components, the reader is referred to \cite{JS-DZ-RO-AS-TS-JR:15,JS:15-thesis}. 

To derive a compact model representation of the power system, it is convenient to introduce the matrices
\begin{equation*}
\begin{split}
&D=\diag(D_i)\in\real^{n\times n},\, M=\diag(M_i)\in\real^{n\times n},\\
\end{split}
\end{equation*}
and the vectors
\begin{equation*}
\begin{split}
&\theta=\col(\theta_i)\in\real^n,\,\omega=\col(\omega_i)\in\real^n,\\
&P^{\text{net}}=\col(P_i^d-G_{ii}V_i^2)\in\real^n,\,
u=\col(u_i)\in\real^n.
\end{split}
\end{equation*}
Also, we introduce the {\em potential function} $U:\real^n \to \real,$
$$U(\theta) = -\sum\nolimits_{\{i,k\}  \in [\mathcal N]^2} |B_{ik}| V_{i}V_{k}\cos(\theta_{ik}).$$ 
Then, the dynamics \eqref{eq:model}, $\forall i\in\mc N,$ can be compactly written as 
\begin{equation}
	\begin{split}
\dot \theta& = \omega,\\
	M \dot \omega & = - D (\omega-\vectorones[n]\omega^d) + P^{\text{net}} +u - \nabla_\theta U(\theta).
\end{split}
\label{eq:model -- 2}%
\end{equation}
Observe that 
due to symmetry of the power flows $P_i,$ 
\begin{equation}
\vectorones[n]^{\top}\nabla_{\theta} U(\theta) = 0.
\label{eq: sumP}
\end{equation}

\section{Nominal DAI-Controlled Power System Model with Fixed Communication Topology and no Delays}
\label{sec: closed-loop model} 
In this section, the employed secondary control scheme is motivated and introduced. Subsequently, we derive the considered resulting {\em nominal} closed-loop system, i.e., without delays and switched topology. For this model, we construct a suitable error system and a strict Lyapunov function in Section~\ref{sec:stabNom}, both of which are instrumental to establish the robust stability results under communication uncertainties in Section~\ref{sec:stabDel}.

\subsection{Secondary Frequency Control: Objectives and Distributed Averaging Integral (DAI) Control}
\label{subsec: sec ctrl}
The whole power system is designed to work at, or at least very close to, the nominal network frequency $\omega^d$ \cite{PK:94}. However, by inspection of a synchronized solution (i.e., a solution with constant uniform frequencies $\omega^*=\omega^s\vectorones[n],$ constant control input $u^*$ and constant phase angle differences $\theta^*_{ik}$) of the system \eqref{eq:model -- 2}, we have that
\begequ
0=\vectorones[n]^\top  M \dot \omega= \vectorones[n]^\top (- D \vectorones[n](\omega^s-\omega^d) + P^{\text{net}} +u^* - \nabla_\theta U(\theta^*)), 
\notag
\endequ
which with \eqref{eq: sumP} implies that
\begequ
\omega^s=\omega^d+ \frac{\vectorones[n]^\top ( P^{\text{net}} +u^*)}{\vectorones[n]^\top D \vectorones[n]}.
\label{omegastar}
\endequ
Note that the loads $G_{ii}V_i^2$ contained in $P^{\text{net}}$ are usually unknown. Hence, in general $\vectorones[n]^\top  P^{\text{net}}\neq0$ and, thus, $\omega^s\neq\omega^d,$ unless  the additional control signal $u^*$ accounts for the power imbalance. Control schemes which yield such $u$ and, consequently, ensure $\omega^s=\omega^d$ are termed {\em secondary frequency controllers}.

Building upon \cite{JWSP-FD-FB:12u,CZ-EM-FD:15,JS-FD:15}, we consider the following secondary frequency control scheme for the system \eqref{eq:model -- 2}
\begin{equation}
\begin{split}
u&=-p,\\
\dot p&=K(\omega-\vectorones[n]\omega^d) - KA\mc L A p.
\end{split}
\label{eq:dai2:0}
\end{equation}
We refer to the control \eqref{eq:dai2:0} as {\em distributed averaging integral (DAI) control law} in the sequel.
Here, $K\in\real^{n\times n}$ is a diagonal gain matrix with positive diagonal entries, $\mc L=\mc L^\top  \in\real^{n\times n}$ is the Laplacian matrix of the undirected and connected communication graph over which the individual generation units can communicate with each other, and $A$ is a positive definite diagonal matrix with the element $A_{ii}>0$ being a coefficient accounting for the cost of secondary control at node $i$. 
It has been shown in \cite{JWSP-FD-FB:12u,CZ-EM-FD:15,JS-FD:15} that the control law \eqref{eq:dai2:0} is a suitable secondary frequency control scheme for the system \eqref{eq:model -- 2}, i.e., it can achieve $\omega^s=\omega^d$ despite unknown (constant) loads; see \cite{ST-MB-CDP:16,NM-CDP:15} for variations of the control \eqref{eq:dai2:0} for dynamic load models. In addition to secondary frequency control, the control law \eqref{eq:dai2:0} can also ensure that the power injections of all generation units satisfy the {\em identical marginal cost} requirement in steady-state, i.e.,
\begin{equation}
	A_{ii}u_{i}^{*} = A_{kk}u_{k}^{*} \mbox{ for all } i\in\mc N,\,k \in \mathcal N \,,
	\label{eq: id marg costs}
\end{equation}
where $A_{ii}$ and $A_{kk}$ are the respective diagonal entries of the matrix $A.$

The {\em nominal} closed-loop system resulting from combining \eqref{eq:model -- 2} with \eqref{eq:dai2:0} is given by
\begin{equation}
	\begin{split}
	\dot \theta &= \omega,\\
M \dot \omega  &= - D (\omega-\omega^d\vectorones[n]) + P^{\text{net}} - \nabla_{\theta} U(\theta)-p,
			\\
	\dot p &= K (\omega-\vectorones[n]\omega^d)  -  KA\mc L A p.
\end{split}
\label{eq: closed loop 1}%
\end{equation}

To formalize our main objective, it is convenient to introduce the notion below.
\begin{definition}[\bf Synchronized motion]
	The system \eqref{eq: closed loop 1} admits a synchronized motion if it has a solution for all $t\geq0$ of the form
	\[
	\theta^{*}(t)=\theta_{0}^{*}+\omega^{*}t,\quad \omega^*=\omega^s\vectorones[n], \quad p^*\in\real^n ,\]
	where $\omega^{s}\in\real$ and $\theta_{0}^{*}\in\real^{n}$ such that\[
	|\theta_{0,i}^{*}-\theta_{0,k}^{*}|<\frac{\pi}{2}\quad\forall i\in\mathcal{N},\;\forall k\in\mathcal{N}_{i}.\]
	\label{def:sm}
\end{definition}

Note that \cite[Lemma 4.2]{JS-FD:15} implies that the system \eqref{eq: closed loop 1} has at most one synchronized motion $\col(\theta^*,\omega^*,p^*)$ (modulo $2\pi$). That motion also satisfies the identical marginal cost requirement \eqref{eq: id marg costs} and is characterized by 
\begin{equation}
p^*=\alpha A^{-1}\vectorones[n], \quad \alpha=\frac{\vectorones[n]^\top  P^{\text{net}}}{\vectorones[n]^\top  A^{-1}\vectorones[n]}
\label{pstar}
\end{equation}
and, hence, $\omega^s=\omega^d,$ see \eqref{omegastar} with $u^*=-p^*.$

\subsection{A Useful Coordinate Transformation}
We introduce a coordinate transformation that is fundamental to establish our robust stability results. 
Recall that $\vectorones[n]^\top  \mc L=0.$ This results in an invariant subspace of the $p$-variables, which makes the construction of a strict Lyapunov function for the system \eqref{eq: closed loop 1} difficult. Therefore, we seek to eliminate this invariant subspace through an appropriate coordinate transformation. 
To this end and inspired by \cite{ROS-RM:04,PL-YJ:08,XW-FD-MJ:15a}, we introduce the variables $\bar p\in\real^{(n-1)}$ and $\zeta\in\real$ via the transformation
\begequ
\begin{split}
\begin{bmatrix} \bar p \\ \zeta \end{bmatrix}&=\mc W^\top K^{-\frac{1}{2}} p,\quad
 \mc W=\begin{bmatrix} W & \frac{1}{\sqrt{\mu}}K^{-\frac{1}{2}}A^{-1}\vectorones[n] \end{bmatrix}, \, 
\end{split}
\label{barp}
\endequ
where $\mu=\| K^{-\frac{1}{2}}A^{-1} \vectorones[n]\|_2^2$ and $W\in\real^{n \times (n-1)}$ has orthonormal columns that are all orthogonal to $K^{-\frac{1}{2}}A^{-1}\vectorones[n]$, i.e., $W^\top  K^{-\frac{1}{2}}A^{-1} \vectorones[n]=\vectorzeros[(n-1)]$. Hence, the transformation matrix $\mc W \in \real^{n \times n}$ is orthogonal, i.e., 
\begequ
\mc W \mc W^\top  
\!=\!WW^\top +\frac{1}{\mu}K^{-\frac{1}{2}}A^{-1}\vectorones[n] \vectorones[n]^\top K^{-\frac{1}{2}}A^{-1}\!=\! I_n.
\label{t_prop}
\endequ 
Accordingly, $\bar p$ is a projection of $p$ on the subspace orthogonal to $K^{-\frac{1}{2}}A^{-1}\vectorones[n]$ scaled by $K^{-\frac{1}{2}}.$ 
Furthermore,
\begequ
\zeta=\frac{1}{\sqrt{\mu}}\vectorones[n]^\top  K^{-1}A^{-1} p
\label{zeta1}
\endequ
can be interpreted as the scaled average secondary control injections of the network. Indeed, from \eqref{eq: closed loop 1} together with the fact that $\vectorones[n]^\top  \mc L=0,$ we have that
\begequ
\dot \zeta = \frac{1}{\sqrt{\mu}}\vectorones[n]^\top  K^{-1}A^{-1}  \dot p=\frac{1}{\sqrt{\mu}}\vectorones[n]^\top A^{-1}(\omega-\vectorones[n]\omega^d),
\notag
\endequ
which by integrating with respect to time and recalling \eqref{zeta1} yields
\begequ
\begin{split}
\zeta&=\frac{1}{\sqrt{\mu}}\vectorones[n]^\top A^{-1}(\theta -\theta_0-\vectorones[n]\omega^dt+K^{-1} p_0)\\
&=\frac{1}{\sqrt{\mu}}\vectorones[n]^\top A^{-1}(\theta -\vectorones[n]\omega^dt) +\bar \zeta_{0},
\end{split}
\label{zeta}
\endequ
where 
\begequ
\bar \zeta_0=\sqrt{\mu}^{-1}\vectorones[n]^\top A^{-1}(K^{-1} p_0 -\theta_{0}).
\label{zeta0}
\endequ
As a consequence, the coordinate $\zeta$ can be expressed by means of $\theta$ and the parameter $\bar \zeta_0.$
Hence, 
\begequ
\begin{split}
p&=K^{\frac{1}{2}}\left(W \bar p + \frac{1}{\sqrt{\mu}} K^{-\frac{1}{2}}A^{-1}\vectorones[n] \zeta \right)\\
&=K^{\frac{1}{2}}\left(W \bar p + \frac{1}{\mu} K^{-\frac{1}{2}}A^{-1}\vectorones[n](\vectorones[n]^\top A^{-1}(\theta-\vectorones[n]\omega^dt)+\bar \zeta_{0}) \right)\,.
\end{split}
\label{pbartop}
\endequ
Accordingly, we define the matrix
\begin{equation*}
\bar {\mc L}=W^\top  K^{\frac{1}{2}}A\mc LAK^{\frac{1}{2}}W\in\real^{(n-1)\times(n-1)},
\end{equation*}
that corresponds to the communication Laplacian matrix $\mc L$ after scaling and projection.
Note that $\mc L$ is connected by assumption, and thus $\bar {\mc L}$ is positive definite. 

In the reduced coordinates, the dynamics \eqref{eq: closed loop 1} become%
	\begequ
\label{eq: closed loop 3}%
\begin{split}
\dot \theta &= \omega,
\\
M \dot {\omega} &= - D (\omega-\omega^d\vectorones[n]) \!+\! P^{\text{net}}\! -\! \nabla_{\theta} U(\theta)\\
&-\! K^{\frac{1}{2}}\!\left(\! W \bar p +\! \frac{1}{\mu} K^{-\frac{1}{2}}A^{-1}\vectorones[n](\vectorones[n]^\top A^{-1}(\theta\!-\!\vectorones[n]\omega^dt)+\bar \zeta_{0})  \!\right)\!,\\
	\dot {\bar p} &= W^\top  K^{-\frac{1}{2}} \dot p\\
	&=W^\top K^{\frac{1}{2}}(\omega-\vectorones[n]\omega^d)- \bar{\mc L}\bar p,
\end{split}%
\end{equation}%
where we have used \eqref{pbartop} and the fact that $ \mc L\vectorones[n]=\vectorzeros[n].$
Note that the transformation \eqref{barp} removes the invariant subspace $\text{span}(\vectorones[n])$ of a synchronized motion of \eqref{eq: closed loop 1} in the $\theta$-variables and shifts it to the invariant subspace
$\text{span}(K^{-\frac{1}{2}}A^{-1})$, where it is factored
out in the orthogonal reduced-order $\bar p$-variables.

\subsection{Error States}
For the subsequent analysis, we make the following standard assumption \cite{JS-RO-AA-JR-TS:14,JS-FD:15}.

\begin{assumption}[\bf Existence of synchronized motion]
The closed-loop system \eqref{eq: closed loop 3} possesses a synchronized motion.
\label{ass:feas} 
\oprocend
\end{assumption}

With initial time $t_0=0,$ as well as
$$
\bar p^* = W^\top K^{-\frac{1}{2}} p^*,
$$
see \eqref{barp}, we introduce the error coordinates
\begequ
\begin{split}
\tilde \omega&=\omega-\vectorones[n]\omega^d,\quad \tilde p=\bar p-\bar p^*,\\
\tilde \theta &= \theta-\theta^*=\theta_0-\theta_0^*+\int_{0}^t \tilde \omega(s)ds,\\
x&=\col\left( \tilde \theta, \, \tilde \omega,\,  \tilde p\right).
\end{split}
\notag
\endequ
Then the dynamics \eqref{eq: closed loop 3} become
\begequ
\begin{split}
\dot{ \tilde \theta} =& {\tilde \omega},\\
M \dot {\tilde \omega} =& - D \tilde \omega - \nabla_{\tilde \theta} U(\tilde \theta+\theta^*)+ \nabla_{\tilde \theta} U(\theta^*)\\
&-\! K^{\frac{1}{2}}W \tilde p -\frac{1}{\mu} A^{-1}\vectorones[n] \vectorones[n]^\top A^{-1}\tilde \theta,\\
	\dot {\tilde p} =& W^\top  K^{\frac{1}{2}} \tilde\omega-  \bar {\mc L} \tilde p
\end{split}%
\label{eq: closed loop 4}%
\end{equation}%
and $x^*=\vectorzeros[(3n-1)]$ is an equilibrium point of \eqref{eq: closed loop 4}. Furthermore, in error coordinates the potential function $U:\real^n \to \real$ reads
$$U(\tilde \theta+\theta^*) = -\sum\nolimits_{\{i,k\}  \in [\mathcal N]^2} |B_{ik}| V_{i}V_{k}\cos(\tilde \theta_{ik}+\theta^*_{ik})$$ 
with
$$
\nabla_{\tilde \theta} U(\tilde \theta+\theta^*)=\frac{\partial U(\tilde \theta+\theta^*)}{\partial \tilde \theta}, \; \nabla_{\tilde \theta} U(\theta^*)=\frac{\partial U(\tilde \theta+\theta^*)}{\partial \tilde \theta}\big|_{\tilde \theta=\vectorzeros[n]}.
$$

Recall from Section~\ref{subsec: sec ctrl} that any synchronized motion of the system \eqref{eq: closed loop 1} satisfies $\omega^s=\omega^d$ and that $p^*$ is uniquely given by \eqref{pstar}. 
Thus, for a fixed value of $\bar \zeta_0$ it follows from \eqref{zeta0} together with \eqref{pbartop} that asymptotic stability of $x^*$ implies convergence of the solutions $\col(\theta,\omega,p)$ of the original system \eqref{eq: closed loop 1} with initial conditions that satisfy
$$
\bar \zeta_0=\sqrt{\mu}^{-1}\vectorones[n]^\top A^{-1}(K^{-1} p_0 -\theta_{0})
$$
to a synchronized motion $\col(\theta^*,\omega^*,p^*),$ the initial angles $\theta_0^*$ of which satisfy 
$$
\bar \zeta_0=\sqrt{\mu}^{-1}\vectorones[n]^\top A^{-1}(K^{-1} p^* -\theta^*_{0}).
$$
As this holds true for any value of $\bar \zeta_0$ and the dynamics \eqref{eq: closed loop 4} are independent of $\bar \zeta_0,$ asymptotic stability of $x^*$ implies convergence of all solutions of the original system \eqref{eq: closed loop 1} to a synchronized motion.

\section{Stability Analysis of the Nominal Closed-Loop System with a Strict Lyapunov Function}
\label{sec:stabNom}
To pave the path for the analysis in Section~\ref{sec:stabDel}, we start by investigating stability of an equilibrium $x^*$ of the {\em nominal} (without delays and with constant communication topology) closed-loop system \eqref{eq: closed loop 4}.

The proposition below provides a stability proof for an equilibrium of the system \eqref{eq: closed loop 4} by employing the following strict Lyapunov function candidate
\begequ
\begin{split}
V=&U(\tilde \theta  + \theta^{*})- \tilde \theta^\top \nabla_{\tilde \theta} U(\theta^*) + \frac{1}{2} \tilde \omega^\top  M \tilde \omega\\
&+ \frac{1}{2\mu} (\vectorones[n]^\top A^{-1}\tilde \theta)^2 + \frac{1}{2} \tilde p^\top  \tilde p\\
&+\epsilon \tilde \omega^\top  AM (\nabla_{\tilde \theta} U(\tilde \theta+\theta^*)-\nabla_{\tilde \theta} U(\theta^*) ),
\end{split}
\label{V}
\endequ
where $\epsilon>0$ is a positive real and sufficiently small parameter. The Lyapunov function \eqref{V} is based on the classic kinetic and potential energy terms  $\omega^{\top} M \omega$ and $U(\theta)$ \cite{MAP:89} written in error coordinates, a Bregman construction to center the Lyapunov function as in \cite{NM-CDP:15,ST-MB-CDP:16}, a Chetaev-type cross term between the (incremental) potential and kinetic energies \cite{bullo2004geometric}, and a quadratic term for the secondary control inputs (also in error coordinates). We have the following result.

\begin{proposition}[\bf Stability of the nominal system]
\label{prop: nom-stab}
Consider the system \eqref{eq: closed loop 4} with Assumption~\ref{ass:feas}.
The function $V$ in \eqref{V} is a strict Lyapunov function for the system \eqref{eq: closed loop 4}.
Furthermore, $x^*=\vectorzeros[(3n-1)]$ is locally asymptotically stable.
\oprocend
\label{stab}
\end{proposition}
\begin{pf}
We first show that the function $V$ in \eqref{V} is locally positive definite. It is easily verified that
$$
\nabla V\big|_{x^*}=\vectorzeros[(3n-1)].
$$
Moreover, the Hessian of $V$ evaluated at $x^*$ is given by
\begin{equation*}
\nabla V^2\big|_{x^*}=\begin{bmatrix}\nabla_{\tilde \theta}^2 U|_{x^*}+\frac{1}{\mu} A^{-1}\vectorones[n]\vectorones[n]^\top A^{-1} & \frac{\epsilon}{2} E_{12}   & \vectorzeros\\  \ast  & M & \vectorzeros\\ \ast &\ast & I_{(n-1)}\end{bmatrix},
\end{equation*}
where
$$
E_{12}=AM\nabla_{\tilde \theta}^2 U|_{x^*}+\nabla_{\tilde \theta}^2 U|_{x^*}MA
$$
and $\vectorzeros$ denotes a zero matrix of appropriate dimension.
Under the standing assumptions, $\nabla_{\tilde \theta}^2 U|_{x^*}$ is a Laplacian matrix of an undirected connected graph. Hence, $\nabla_{\tilde \theta}^2 U|_{x^*}$ is positive semidefinite with $\ker(\nabla_{\tilde \theta}^2 U|_{x^*})=\text{span}(\vectorones[n]).$ Furthermore, $A^{-1}\vectorones[n] \vectorones[n]^\top A^{-1}$ is positive semidefinite and $\ker(A^{-1}\vectorones[n] \vectorones[n]^\top A^{-1})\cap \ker(\nabla_{\tilde \theta}^2 U|_{x^*})=\vectorzeros[n].$ In addition, $M$ is a diagonal matrix with positive diagonal entries. Thus, all block-diagonal entries of $\nabla V^2\big|_{x^*}$ are positive definite. This implies that there is a sufficiently small $\epsilon^{*}>0$ such that for all $\epsilon \in {]0,\epsilon^{*}]}$ we have that $\nabla V^2\big|_{x^*}>0.$ We choose such $\epsilon.$ Therefore, $x^*$ is a strict minimum of $V.$

The time derivative of $V$ along solutions of \eqref{eq: closed loop 4} is given by
\begequ
\begin{split}
\!\!\!\dot {V}
\!&\!=-\xi^\top \! \! \begin{bmatrix} \epsilon A & 0.5 \epsilon AD & 0.5 \epsilon AK^{\frac{1}{2}}W \!\\
	\ast& D\!-\!0.5\epsilon E_{22} & \vectorzeros[n\times (n-1)]\\
	\ast  & \ast & \bar{\mc L}
\end{bmatrix}\!
\xi,\!\!
\end{split}
\label{dotmcV1}
\endequ
where we have used the property that
$$
(\nabla_{\tilde \theta} U(\tilde \theta+\theta^*)-\nabla_{\tilde \theta} U(\theta^*) )^\top  \vectorones[n]=0
$$
and defined the shorthand
\begequ
E_{22}=(AM\nabla_{\tilde \theta}^2 U(\tilde \theta+\theta^*)+\nabla_{\tilde \theta}^2 U(\tilde \theta+\theta^*) AM),
\label{e22}
\endequ
as well as
$$
\xi=\col\left(\nabla_{\tilde \theta} U(\tilde \theta+\theta^*)-\nabla_{\tilde \theta} U(\theta^*),\, \tilde \omega,\, \tilde p\right).
$$
Note that the matrix $E_{22}$ only depends on the cosines of the angles $\tilde \theta.$ Hence, there is a positive real constant $\gamma,$ such that
$
E_{22} \leq \gamma I_n$ for all $\tilde \theta \in\real^n.
$
Furthermore, $A,$ $D$ and $\bar {\mc L}$ are positive definite matrices.
Thus, Lemma \ref{Lemma: matrix regularization} in Appendix~\ref{appendix2} implies that there is a sufficiently small $\epsilon^{**} \in {(0,\epsilon^{*}]}$ such that for any $\epsilon \in {(0,\epsilon^{**}]}$ the matrix on the right-hand side of \eqref{dotmcV1} is positive definite and, thus, $\dot{ V}<0$ for all $\xi\neq \vectorzeros[(3n-1)].$ 
Consequently, by choosing such $\epsilon,$ $V$ is a strict Lyapunov function for the system \eqref{eq: closed loop 4} and, by Lyapunov's theorem \cite{HKK:02}, $x^*$ is asymptotically stable, completing the proof.
\oprocend
\end{pf}

\section{Stability of the Closed-Loop System with Dynamic Communication Topology and Delays}
\label{sec:stabDel}

This section is dedicated to the analysis of cyber-physical aspects in the form of communication uncertainties on the performance of the DAI-controlled power system model \eqref{eq: closed loop 4}. 

\subsection{Modeling and Problem Formulation}
\label{subsec: prob stat}
Recall from Section~\ref{Section: Introduction} that the most relevant practical communication uncertainties in the context of DAI control are message delays, message losses and link failures. We follow the approaches in \cite{ROS-RM:04,ROS-AF-RM:07,PL-YJ:08,EF:14,EF:14b} to model these phenomena.

Thus, link failures and packet losses are modeled by a dynamic communication network with switched communication topology $\mc G_{\sigma(t)},$ where $\sigma:\real_{\geq0}\to \mc M$ is a switching signal and $\mc M=\{1,2,\ldots,\nu\},$ $\nu\in\real_{>0},$ is an index set.
The {\em finite} set of all possible network topologies amongst $|\mc N|=n$ nodes is denoted by 
$
\Gamma=\{\mc G_1, \mc G_2,\ldots \mc G_\nu \}.
$
The Laplacian matrix corresponding to the  index $\ell=\sigma(t)\in\mc M$
is denoted by
$\mc L_\ell = \mc L_{\ell}^{\top}=\mc L(\mc G_\ell)\in\real^{n\times n}.$ 
We employ the following standard assumption on $\mc G_{\sigma(t)}$ \cite{ROS-RM:04,ROS-AF-RM:07}.
\begin{assumption}{\bf (Uniformly connected communication topologies)}
The communication topology $\mc G_{\sigma(t)}$ is undirected and connected for all $t\in\real_{\geq0}.$
\label{Ass: Communication topology}
\oprocend 
\end{assumption}
With regards to communication delays, we suppose that a message sent by generation unit $k\in \mc N$ to the generation unit $i\in\mc N$ over the communication channel (i.e., edge) $\{i,k\}$ is affected by a fast-varying delay 
$\tau_{ik}:\real_{\geq0}\to[0,h],$ $h\in\real_{\geq0}$, where the qualifier "fast-varying" means that there are no restrictions imposed on the existence, continuity, or boundedness of $\dot \tau_{ik}(t)$  \cite{EF:14,EF:14b}.
The resulting control error $e_{ik}$ is then computed as
$$
e_{ik}(t)=A_{ii}p_i(t-\tau_{ik}(t))-A_{kk}p_{k}(t-\tau_{ik}(t)),
$$
i.e., the protocol is only executed after the message from node $k$ arrives at node $i.$ Note that we allow for asymmetric delays, i.e., $\tau_{ik}(t)\neq \tau_{ki}(t).$ Furthermore, as standard in sampled-data networked control systems \cite{EF:14,EF:14b}, the delay $\tau_{ik}(t)$ may be piecewise-continuous in $t$. Also, we assume that the switches in topology don't modify the delays between two connected nodes.

In order to write the resulting closed-loop system compactly, we introduce the matrices $T_{\ell,m} \in \real^{n\times n},$ $m=1,\ldots,2|\mc E_\ell|,$ where $|\mc E_\ell|$ is the number of edges of the undirected graph with index $\ell=\sigma(t)\in\mc M$ of the dynamic communication network of the DAI control \eqref{eq:dai2:0}, $m$ denotes the information flow from node $i$ to $k$ over the edge $\{i,k\}$ with delay $\tau_m=\tau_{ik}$ and all elements of $T_{\ell,m}$ are zero besides the entries
\begequ
t_{\ell,m,ii}=1,\quad t_{\ell,m,ik}=-1.
\label{tm}
\endequ
As we allow for $\tau_{ik}\neq \tau_{ki},$ we require $2|\mc E_\ell|$ matrices $T_{\ell,m}$ in order to distinguish between the delayed information flow from $k$ to $i$ (with delay $\tau_{ik}$) and that from $k$ to $i$ (with delay $\tau_{ki}$). Note that by summing over all $T_{\ell,m}$ we recover the full Laplacian matrix of the communication network corresponding to the topology index $\ell=\sigma(t)\in\mc M$, i.e.,
\begequ
\mc L_\ell=\sum_{m=1}^{2|\mc E_\ell|} T_{\ell,m}.
\notag
\endequ

With the above considerations, the closed-loop system \eqref{eq: closed loop 1} becomes the switched nonlinear delay-differential system 
\begequ
\begin{split}
	\dot \theta &= \omega,\\
	M \dot \omega  &= - D (\omega-\omega^d\vectorones[n]) + P^{\text{net}} - \nabla_{\theta} U(\theta)-p,\\
	\dot p &= K (\omega-\vectorones[n]\omega^d)  -  KA \left(\sum_{m=1}^{2|\mc E_\ell|} T_{\ell,m}A p(t-\tau_m)\right),
\end{split}%
\label{eq: closed loop 5}%
\end{equation}%
where $\tau_m(t)\in[0,h_m],$ $m=1,\ldots,2|\mc E_\ell|$ are fast-varying delays and $T_{\ell,m}$ corresponds to the $m$-th delayed (directed) channel of the $\ell$-th communication topology corresponding to the topology index $\ell=\sigma(t)\in\mc M$ of the dynamic communication network of the DAI control \eqref{eq:dai2:0}. 

We are interested in the following problem.
\newpage 

\begin{problem}{(\bf Conditions for robust stability)}
\label{prob: robust stability}
Consider the system \eqref{eq: closed loop 5}. Given $h_m\in\real_{\geq0},$ $m=1,\ldots,2\bar {\mc E},$ $\bar {\mc E}=max_{\ell\in\mc M}{|\mc E_\ell|},$ derive conditions under which the solutions of the system \eqref{eq: closed loop 5} converge asymptotically to a synchronized motion.
\oprocend
\label{problem}
\end{problem}
As in the nominal scenario, we make the assumption below.%
\begin{assumption}[\bf Existence of synchronized motion]
The closed-loop system \eqref{eq: closed loop 5} possesses a synchronized motion.
\label{ass:feas2} 
\oprocend
\end{assumption}

From \eqref{pbartop} it follows that for any $m=1,\ldots,2\bar {\mc E},$
$$
p(t-\tau_m)=\!K^{\frac{1}{2}}\!\left(\!W \bar p(t-\tau_m) + \frac{1}{\sqrt{\mu}} K^{-\frac{1}{2}}A^{-1}\vectorones[n] \zeta(t-\tau_m)\! \right),
$$
which together with the fact that by construction, see \eqref{tm}, $T_{\ell,m}\vectorones[n]=\vectorzeros[n]$ implies that
$$
T_{\ell,m} A p(t-\tau_m)=T_{\ell,m}AK^{\frac{1}{2}}W \bar p(t-\tau_m).
$$
Hence, with Assumption~\ref{ass:feas2} and by following the steps in Section~\ref{sec: closed-loop model}, we represent the system \eqref{eq: closed loop 5} in reduced-order error coordinates as, cf., \eqref{eq: closed loop 4},
\begequ
\begin{split}
\dot{ \tilde \theta} =& {\tilde \omega},\\
M \dot {\tilde \omega} =& - D \tilde \omega - \nabla_{\tilde \theta} U(\tilde \theta+\theta^*)+ \nabla_{\tilde \theta} U(\theta^*)\\
&-\! K^{\frac{1}{2}}W \tilde p -\frac{1}{\mu} A^{-1}\vectorones[n] \vectorones[n]^\top A^{-1}\tilde \theta,\\
	\dot {\tilde p} =& W^\top  K^{\frac{1}{2}} \tilde\omega- \sum_{m=1}^{2|\mc E|} \mc T_{\ell,m} \tilde p(t-\tau_m),
\end{split}%
\label{eq: closed loop 6}%
\end{equation}%
where we defined, analogous to $\bar {\mc L}$,
\begequ
\mc T_{\ell,m} = W^\top  K^{\frac{1}{2}}AT_{\ell,m} AK^{\frac{1}{2}}W \in\real^{(n-1)\times(n-1)},
\label{mcTell}
\endequ 
which satisfies
\begequ
\sum_{m=1}^{2|\mc E_\ell|}\mc T_{\ell,m} = W^\top  K^{\frac{1}{2}}A\mc L_{\ell} AK^{\frac{1}{2}}W=:\bar{\mc L}_\ell.
\label{barLell}
\endequ
Note that, by assumption, the graph associated to $\mc L_\ell$ is connected and thus $\bar {\mc L_\ell}$ is positive definite for any $\ell\in\mc M$.

By construction, the system \eqref{eq: closed loop 6} has an equilibrium point $z^*=\col(\tilde\theta,\tilde\omega,\tilde p)=\vectorzeros[(3n-1)].$ Furthermore, by the same arguments as in Section~\ref{subsec: sec ctrl},
asymptotic stability of $z^*$ implies convergence of the solutions of the system \eqref{eq: closed loop 5} to a synchronized motion.
As a consequence of this fact, we provide a solution to Problem~\ref{problem} by studying stability of $z^*.$

\subsection{Stability of the Closed-Loop System with Dynamic Communication Topology and Delays}
\label{sec:stabConst}

We analyze stability of equilibria of the system \eqref{eq: closed loop 6} with dynamic communication topology and time delays. The result is formulated for fast-varying piecewise-continuous bounded delays $\tau_m(t)\in[0,h_m]$ with $h_m\in\real_{\geq 0},$ $m=1,\ldots,2\bar {\mc E}.$ 
For the case of constant delays, stability can be verified via the same conditions.

To streamline our main result, we recall from Section~\ref{sec: closed-loop model} that the matrix $A$ can be used to achieve the objective of identical marginal costs. 
As a consequence, the choice of $A$ influences the corresponding equilibria of the system \eqref{eq: closed loop 6}, see \eqref{pstar}. But the equilibria are independent of the integral control gain matrix $K.$
Hence, we may chose $K$ as a free tuning parameter to ensure stability of the closed-loop system \eqref{eq: closed loop 6}. The proof of the proposition below is given in Appendix~\ref{appendix:constant}.
\begin{proposition}[\bf Robust stability]
	Consider the system \eqref{eq: closed loop 6} with Assumptions~\ref{Ass: Communication topology} and \ref{ass:feas2}. Fix $A$ and $D$ as well as some $h_m\in\real_{\geq 0},$ $m=1,\ldots,2\bar {\mc E}.$  Select $K$ such that for all $\mc T_{\ell,m}$ and $\bar{\mc L}_{\ell}$ defined in \eqref{mcTell}, respectively \eqref{barLell}, $\ell=1,\ldots,|\mc M|,$ there exist matrices $S_m>0\in\real^{(n-1)\times(n-1)},$ $R_m>0\in\real^{(n-1)\times (n-1)}$ and $S_{12,m}\in\real^{(n-1)\times (n-1)}$ satisfying
\begequ
\Psi=\begin{bmatrix}
	\Psi_{11} &  \Psi_{12} & \vectorzeros & \Psi_{14}\\
	\ast  &  \Psi_{22} & \Psi_{23} &  \Psi_{24}\\
	\ast  & \ast  & R+S & S_{12}+S\\                                                
	\ast  & \ast  &  \ast & R+S+ \bar \Psi_{44}
\end{bmatrix}>0,
\label{psi}
\endequ
where
\begequ
\begin{split}
	R&=\text{blockdiag}(R_m),\; S=\text{blockdiag}(S_m),\\
	S_{12}&=\text{blockdiag}(S_{12,m}),\; \bar R=\sum_{j=1}^{2\bar {\mc E}} h_j^2 R_j,\\
	\Psi_{11}&=D - K^{\frac{1}{2}}W \bar R W^\top K^{\frac{1}{2}},\;
	\Psi_{22}=\bar{\mc L}_{\ell}- \bar{\mc L}_{\ell}\bar R\bar{\mc L}_{\ell},\\
	\bar \Psi_{44}&=\text{blockdiag}\left(-\mc T_{\ell,m}^\top \bar R \mc T_{\ell,m} \right),\\
	\Psi_{12}&=K^{\frac{1}{2}}W \bar R\bar{\mc L}_{\ell},\;
	\Psi_{23}=-\begin{bmatrix} S_1 & \ldots &  S_{2\bar {\mc E}}\end{bmatrix},\\
	\Psi_{14}&=\begin{bmatrix} \bar \Psi_{14,1} & \ldots &  \bar \Psi_{14,{2\bar {\mc E}}}\end{bmatrix},\;\Psi_{24}=\begin{bmatrix} \bar \Psi_{24,1} & \ldots &  \bar \Psi_{24,{2\bar {\mc E}}}\end{bmatrix},\\
	\bar \Psi_{14,m}&=-K^{\frac{1}{2}}W \bar R\mc T_{\ell,m},\\
	\bar \Psi_{24,m}&=\bar{\mc L}_{\ell} \bar R \mc T_{\ell,m}-S_m-0.5\mc T_{\ell,m},
\end{split}
	\label{psiii}
\endequ 
 $\vectorzeros$ denotes a zero matrix of appropriate dimensions and
 \begequ
 \begin{bmatrix}
 	R & S_{12}\\ \ast & R
 \end{bmatrix}\geq0.
 \label{rs12}
 \endequ
Then the equilibrium $z^*=\vectorzeros[(3n-1)]$ is locally uniformly asymptotically stable for all fast-varying delays $\tau_m(t)\in[0, h_m].$
\label{stab2}
	\oprocend
\end{proposition}
The stability certificate \eqref{psi}, \eqref{rs12} is based on a LKF derived from the Lyapunov function \eqref{V}, and it is fairly tight (see Section~\ref{sec: sim}). Note that the evaluation of the certificate \eqref{psi}, \eqref{rs12} and the corresponding controller tuning inherently is centralized and requires complete system information. 
However, the stability certificate \eqref{psi}, \eqref{rs12} can be also made wieldy for a practical plug-and-play control implementation by trading off controller performance for robust stability. The  following corollary makes this idea precise for the case of uniform delays, i.e., $\tau_m(t)=\tau(t)\in [0,h],$ $h_m=h,$ $h\in\real_{\geq0},$ $m=1,\ldots,2\bar{\mc E}.$

\begin{corollary}[Performance-robustness-trade-off]
	Consider the system \eqref{eq: closed loop 6} with Assumptions~\ref{Ass: Communication topology} and \ref{ass:feas2}. Fix $A$ and $D$ as well as some $h\in\real_{\geq 0}.$  Suppose that $\tau_m(t)=\tau(t)\in [0,h]$ and that \eqref{rs12} is satisfied with strict inequality. Set $K=\kappa \mc K,$ where $\kappa\in\real_{\geq0}$ and $\mc K\in\real_{>0}^{n\times n}$ is a diagonal matrix with positive diagonal entries. Then there is $\kappa>0$ sufficiently small, such that the equilibrium $z^*=\vectorzeros[(3n-1)]$ is locally uniformly asymptotically stable for all fast-varying delays $\tau(t)\in[0, h].$
	\label{cor:stab}
\end{corollary}
\begin{pf}
With $K=\kappa \mc K,$ we have that 
$$
\mc T_{\ell,m} = \kappa \hat{\mc T}_{\ell,m},\quad \bar{\mc L}_{\ell}=\kappa \hat{\mc L}_{\ell}=\kappa\sum_{m=1}^{2\bar{\mc E}} \hat{\mc T}_{\ell,m}, 
$$
with $\mc T_{\ell,m}$ and $\bar{\mc L}_{\ell}$ defined in \eqref{mcTell}, respectively \eqref{barLell}.  Furthermore, we write the free parameter matrix $S$ as $S=\kappa \mc S,$ $\mc S>0.$ 
Consequently, the matrices in \eqref{psiii} also become $\kappa$-dependent. Then, for $\tau_m(t)=\tau(t)\leq h$ the matrix $\Psi$ in \eqref{psi} can be written as
 \begequ
\Psi=\begin{bmatrix}
	D& \vectorzeros & \vectorzeros & \vectorzeros\\
	\ast & \vectorzeros  & \vectorzeros &  \vectorzeros\\
	\ast & \ast & R & S_{12}\\                                                
	\ast & \ast &  \ast & R
\end{bmatrix} + \kappa \begin{bmatrix}
\Phi_{11} & \Phi_{12} & \vectorzeros & \Phi_{14}\\
\ast & \Phi_{22}  & -\mc S &   \Phi_{24}\\
\ast & \ast & \mc S & \mc S\\                                                
\ast & \ast &  \ast & \Phi_{44}
\end{bmatrix},
\label{phi}
\endequ 
where
\begequ
\begin{split}
\Phi_{11}&= -h^2\mc K^{\frac{1}{2}}WRW^\top \mc K^{\frac{1}{2}},\, \Phi_{12}=h^2\kappa^{\frac{1}{2}}\mc K^{\frac{1}{2}}W R\hat{\mc L}_{\ell} ,\\
\Phi_{22}&=\hat{\mc L}_{\ell}- h^2\kappa\hat{\mc L}_{\ell} R \hat{\mc L}_{\ell},\, \Phi_{44}=\mc S-h^2\kappa \bar{\mc L}_{\ell} R \bar{\mc L}_{\ell},\\
\Phi_{14}&=-h^2\kappa^{\frac{1}{2}} \mc K^{\frac{1}{2}}W R\hat{\mc L}_{\ell},\,\Phi_{24}=h^2\kappa \hat{\mc L}_{\ell} R \hat{\mc L}_{\ell}-\mc S-0.5\hat{\mc L}_{\ell}.
\end{split}
\label{phi_entries}
\endequ
Note that $\Phi_{22}$ can be written as
\begequ
\hat{\mc L}_{\ell}\left(\hat{\mc L}_{\ell}^{-1}- h^2\kappa R \right) \hat{\mc L}_{\ell}.
\label{phi22}
\endequ
Hence, by continuity, for any given $\hat{\mc L}_{\ell}>0$ and $R>0$ we can find a small enough $\kappa,$ such that $\Phi_{22}>0.$
In addition, $D>0$ and \eqref{rs12} is satisfied with strict inequality  by assumption. Therefore, the matrix $\Psi$ in \eqref{phi} is a parameter-dependent composite matrix of the form stated in Lemma~\ref{Lemma: matrix regularization}.
Consequently,  Lemma~\ref{Lemma: matrix regularization} implies that for given $h,$ $A,$ $D$ and $\mc L_\ell,$ $\ell=1,\ldots,\nu,$ there is always a sufficiently small gain $\kappa$ such that there exist matrices $\mc S,$ $R$ and $S_{12}$ satisfying conditions \eqref{psi}, \eqref{rs12}. 
\oprocend
\end{pf}

The claim in Corollary~\ref{cor:stab} is in a very similar spirit to the result obtained in \cite{ROS-RM:04,ROS-AF-RM:07} for the standard linear consensus protocol with delays. In essence, Corollary~\ref{cor:stab} shows that there is a trade-off between delay-robustness, i.e., feasibility of conditions \eqref{psi}, \eqref{rs12} and a high gain matrix $K$ for the DAI controller \eqref{eq:dai2:0}. 
Aside from displaying an inherent performance-robustness trade-off, Corollary~\ref{cor:stab} allows us to certify robust stability (for any switched communication topology) based merely on sufficiently small control gains and without evaluating a linear matrix inequality in a centralized fashion.

We remark that it is also possible to obtain a completely decentralized (though in general more conservative) tuning criterion by further bounding the matrices in \eqref{phi}, respectively \eqref{psi}.

\begin{remark}
Conditions \eqref{psi}, \eqref{rs12} are linear matrix inequalities that can be efficiently solved via standard numerical tools, e.g., \cite{JL:04}. 
Furthermore, instead of checking \eqref{psi}, \eqref{rs12} for all $\bar {\mc L}_{\ell}$ it also suffices to do so for their convex hull \cite{QGW:91}. In addition, checking the convex hull gives a robustness criterion for all unknown topology configurations within the chosen convex hull.
Compared to the related results on stability of delayed port-Hamiltonian systems with application to microgrids with delays \cite{JS-EF-RO:15,JS-EF-RO-JS:16}, the conditions \eqref{psi}, \eqref{rs12} are independent of the specific equilibrium point $z^*.$
\oprocend
\end{remark}

\begin{remark}
\label{rem1}
Note that
$
\hat{\mc L}_{\ell}^{-1} \geq {1}/{\lambda_{\max}(\hat{\mc L}_{\ell})}I_{(n-1)},
$
where $\lambda_{\max}(\hat{\mc L}_{\ell})$ denotes the maximum eigenvalue of the matrix $\hat{\mc L}_{\ell},$ $\ell=1,\ldots,\nu.$ Hence, \eqref{phi22} shows that the admissible largest gain $\kappa$ in order for $\Phi_{22}$ (defined in \eqref{phi_entries}) to be positive definite is constrained by the maximum eigenvalue of all matrices $\hat{\mc L}_{\ell}.$ 
This fact is confirmed by the numerical experiments in Section~\ref{sec: sim}. Similar observations have also been made numerically for consensus systems \cite{PL-YJ:08}.
\oprocend 
\end{remark}

\section{Numerical Example}
\label{sec: sim}

\begin{figure}
	\centering
		\includegraphics[width=.475\textwidth]{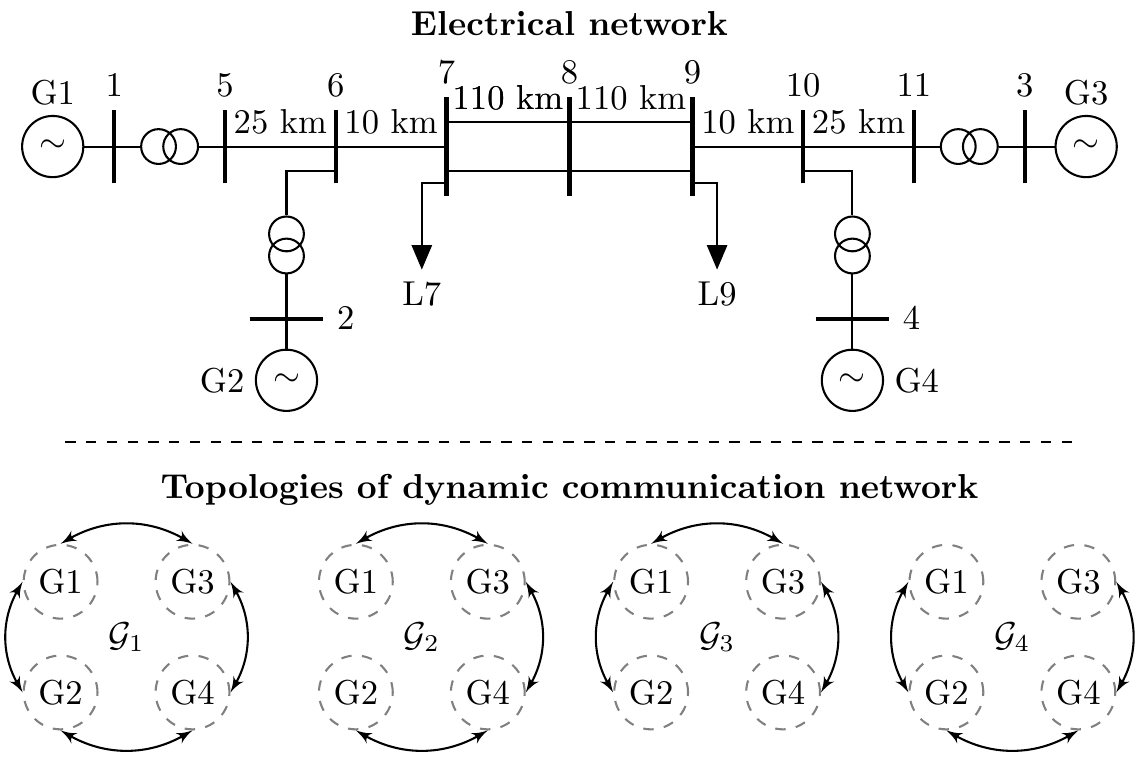}
	\caption{Kundur's two-area-four-machine test system taken from \cite[Example 12.6]{PK:94} and below the four different topologies of the switched communication network.}
	\label{fig:kundur_2area}
\end{figure}

We illustrate the effectiveness of the proposed approach on Kundur's four-machine-two-area test system \cite[Example 12.6]{PK:94}. 
The considered test system is shown in Fig.~\ref{fig:kundur_2area}. The employed parameters are as given in \cite[Example 12.6]{PK:94} with the only difference that we set the damping constants to $D_i=1/(0.05\cdot2\pi\cdot 60)$ pu (with respect to the machine loading \mbox{$S_{\text{SG}}=[700,700,719,700]$} MVA). The system base power is $S_\text{B}=900$ MVA, the base voltage is $V_\text{B}=230$ kV and the base frequency is \mbox{$\omega_\text{B}=1$ rad/s.} 

For the tests, we consider the four different communication network topologies shown in Fig.~\ref{fig:kundur_2area}. Topology $\mc G_1$ is the nominal topology, representing a ring graph. In topology $\mc G_2,$ respectively $\mc G_3$ and $\mc G_4,$ one of the links is broken. Note that each of the configurations is connected and, hence, the dynamic communication network given by the topologies $\Gamma=\{\mc G_1, \mc G_2, \mc G_3, \mc G_4 \}$ satisfies Assumption~\ref{Ass: Communication topology}. Furthermore, we set $A=\diag(S_{\text{SG}}/S_\text{B})$ and $K=\kappa \mc K,$ where $\mc K=0.05A^{-1}$ and $\kappa$ is a free tuning parameter.

We consider an exemplary scenario with fast-varying delays $\tau_m\in[0,h_m]$s, $h_m= 2$s, $m=1,\ldots,2\bar{\mc E}.$ With the chosen parameters, we check conditions \eqref{psi}, \eqref{rs12}. The numerical implementation is done in Yalmip \cite{JL:04}. 
In order to identify the maximum admissible gain $\kappa,$ we select an initial value for $\kappa$ of $\kappa_{\text{init}}=2.0$ and iteratively decrease the value of $\kappa$ until conditions \eqref{psi}, \eqref{rs12} are satisfied.
The obtained feasible gain is $\kappa_{\text{feas}}=1.544.$ 

\begin{figure}
	\centering
		\includegraphics[width=.475\textwidth]{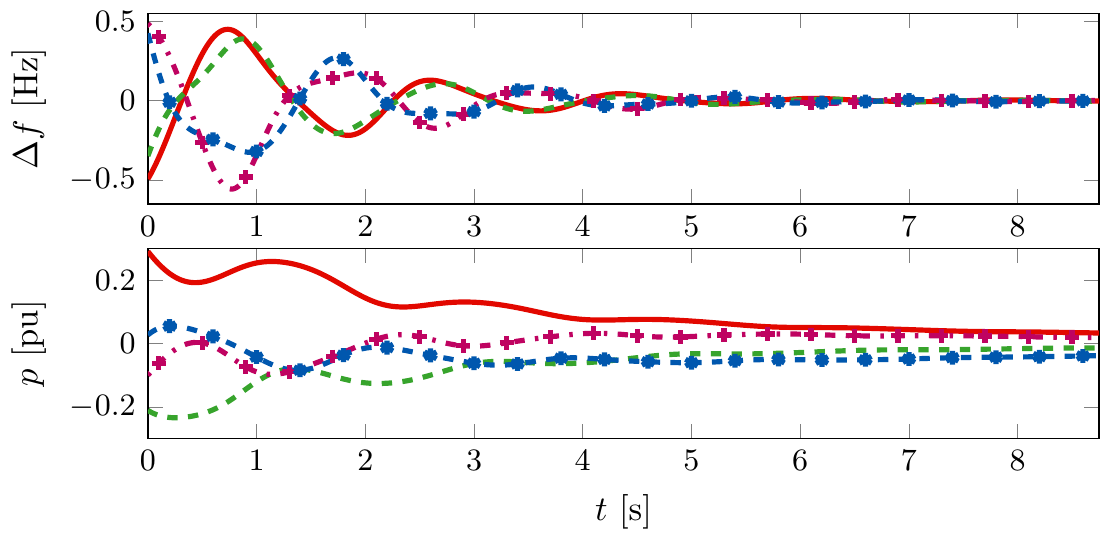}
	\caption{Simulation example for $\kappa=1.544,$ $h_m=2$s, $m=1,\ldots,4$ and arbitrary initial conditions in a neighborhood of $z^{*,1}$. The lines correspond to the following units:
G1 '\textcolor[rgb]{0.88627,0.03137,0.00000}{--}', G2 '\textcolor[rgb]{0.21569,0.64314,0.17255}{- -}', G3 '\textcolor[rgb]{0.74902,0.01176,0.38039}{+-}' and G4 '\textcolor[rgb]{0.00000,0.34118,0.68235}{* -}'. }
	\label{fig:sim}
\end{figure}

The simulation results shown in Fig.~\ref{fig:sim} confirm the fact that the systems' trajectories converge to a synchronized motion if $\kappa_{\text{feas}}=1.544.$ 
Following standard practice in sampled-data networked control systems \cite{EF:14,EF:14b}, the time-varying delays are implemented as piecewise-continuous signals. For the present simulations, we have used the rate transition and variable time delay blocks in Matlab/Simulink with a sampling time $T_s=2$ms to generate the fast-varying delays.
During the simulation, the communication topology switches randomly every $0.5$s. Variations in the switching interval did not lead to meaningful changes in the system's behavior. This shows that the DAI control is very robust to dynamic changes in the communication topology. Furthermore, our experiments confirm the observation in Remark~\ref{rem1} and \cite{PL-YJ:08} that for a given upper bound on the delay the feasibility of conditions \eqref{psi}, \eqref{rs12} is highly dependent on the largest eigenvalue of the Laplacian matrices. 

We verify the conservativeness of the sufficient conditions \eqref{psi}, \eqref{rs12} in simulation. Recall that the conditions are equilibrium-independent. Hence, if feasible they guarantee (local) stability of any equilibrium point satisfying the standard requirement of the equilibrium angle differences being contained in an arc of length $\pi/2,$ cf. Assumption~\ref{ass:feas2}. Therefore, we consider three different operating conditions: at first, the nominal equilibrium point $z^{*,1}$ reported in \cite[Example 12.6]{PK:94} and subsequently two operating points $z^{*,2}$ and $z^{*,3}$ under more stressed conditions, i.e., with some angle differences being closer to the ends of the $\pi/2$-arc. The angles for all three scenarios are given in Table~\ref{tab:xs}.

For the nominal operating point $z^{*,1}$ the maximum feasible gain obtained for fast-varying delays with $h_m=2$s via simulation experiments is $\kappa_{\text{feas,sim}}=6.330.$ For higher values of $\kappa,$ the system exhibits limit cycling behavior. The value of $\kappa_{\text{feas,sim}}=6.330$ is about $4.1$ times larger than the value of $\kappa_{\text{feas}}=1.544$ obtained via Proposition~\ref{stab2}. 
In the case of $z^{*,2}$ and with the same delays the maximum feasible gain identified via simulation is 
$\kappa_{\text{feas,sim}}=2.856,$ which is about $1.85$ times larger as the value of $\kappa_{\text{feas}}=1.544$ obtained from conditions \eqref{psi}, \eqref{rs12}. The maximum feasible gain is further reduced in the last operating scenario with equilibrium $z^{*,3},$ where we obtain $\kappa_{\text{feas,sim}}=1.698.$ 
This value is only $1.1$ times larger than the value of $\kappa_{\text{feas}}=1.544$ obtained from conditions \eqref{psi}, \eqref{rs12}. 
We observe a very similar behavior in the case of constant delays. 
This shows that - for the investigated scenarios - our (equilibrium-independent) conditions are fairly conservative at equilibria corresponding to less stressed operating conditions, while they are almost exact under highly stressed operating conditions. 

\begin{center}
	\begin{table}
		\caption{Values of phase angles at considered equilibria $z^{*,1}$ to $z^{*,3}.$}
		\label{tab:xs}
		\small
		\begin{tabular}{p{1.5cm} | p{6cm}}
			\toprule
			$\theta^{*,1}$ [rad]& $[ 0.224 ,   0.117,   -0.076,   -0.189]\cdot\pi/2$ \\
			$\theta^{*,2}$ [rad]& $[0.3,-0.4,-0.5,0.4]\cdot3\pi/8$ \\
			$\theta^{*,3}$ [rad]& $[0.3,-0.4,-0.5,0.4]\cdot\pi/2$ \\
 			\bottomrule
		\end{tabular}
\end{table}
\end{center}

\section{Conclusions} 
\label{sec:Conclusions}
We have considered the problem of robust stability of a DAI-controlled power system with respect to cyber-physical uncertainties in the form of constant and fast-varying communication delays as well as link failures and packet losses. The phenomena of link failures and packet losses lead to a time-varying communication topology with arbitrary switching. 
For this setup, we have derived sufficient delay-dependent stability conditions by constructing a suitable common LKF. The stability conditions can be verified without knowledge of the operating point and reflect a performance and robustness trade-off, which is in a very similar manner to that of the standard linear consensus protocol with delays investigated in \cite{ROS-RM:04,ROS-AF-RM:07}.

In addition, 
the approach has been applied to Kundur's two-area-four-machine test system.
The numerical experiments show that our sufficient equilibrium-independent conditions are very tight for stressed operating points, i.e., with (some) stationary angle differences close to $\pm\pi/2,$ while they are more conservative for less stressed operating points.

In future research we will seek to relax some of the modeling assumptions made in the paper, e.g., on constant voltage amplitudes, constant load models, and investigate further applications of DAI and related distributed control methods to power systems and microgrids. Another interesting, yet technically challenging, open question is to weaken the connectivity assumption imposed on the DAI communication network, e.g., by using the notion of joint connectivity \cite{AJ-JL-ASM:02}.

\bibliographystyle{elsarticle-num}
\bibliography{./NewRefs,./MainRefs,./FD,./NewRefs1}

\appendix 
{\bf Appendix}
\section{A Matrix Regularization Lemma}
\label{appendix2}
\begin{lem}[\bf Matrix regularization]
\label{Lemma: matrix regularization}
Consider the two symmetric real block matrices
\begin{equation*}
A=\begin{bmatrix}
\vectorzeros[p\times p] & \vectorzeros[p\times q]
\\
\ast  & A_{22}
\end{bmatrix}
\quad \text{ and }\quad 
B = \begin{bmatrix}
B_{11}& B_{12}
\\
\ast & B_{22}
\end{bmatrix},
\end{equation*}
with $A_{22}\in\real^{q\times q},B_{11}\in\real^{p\times p},\, B_{22}\in\real^{q\times q},\,B_{12}\in\real^{p\times q}.$
Suppose that $A_{22}$ and $B_{11}$ are positive definite. 
Then there exists a (sufficiently small) positive real $\epsilon$ such that the composite matrix 
$C_{\epsilon} = A+\epsilon  B
$
is positive definite.
\end{lem}
\begin{pf}
The composite matrix reads as
\begin{equation*}
C_{\epsilon} =
\begin{bmatrix}
\epsilon B_{11}& \epsilon B_{12}
\\
\ast & A_{22}+\epsilon B_{22}
\end{bmatrix}
\,.
\end{equation*}
As $B_{11}$ is positive definite by assumption, $\epsilon B_{11}>0$ for any $\epsilon>0.$
Furthermore, by applying the Schur complement to $C_{\epsilon}$ we obtain 
$$A_{22}+\epsilon\left(B_{22}-B_{12}^\top B_{11}^{-1} B_{12}\right),$$ 
which, as $A_{22}>0$ by assumption, is positive definite for small enough $\epsilon$. The latter together with $\epsilon B_{11}>0$ implies that $C_{\epsilon}>0$ for small enough $\epsilon$, completing the proof. \oprocend
\end{pf}

\section{Proof of Proposition~\ref{stab2}}
\label{appendix:constant}
We give the proof of Proposition~\ref{stab2}.
	The proof is inspired by \cite{PL-YJ:08,EF:14,JS-EF-RO:15,JS-EF-RO-JS:16} and established by constructing a common LKF for  the system \eqref{eq: closed loop 6}. 
	Consider the LKF with $m=1,\ldots,2\bar {\mc E},$
\begequ
\begin{split}
	\mc V=&V + \sum_{m=1}^{2\bar {\mc E}}\mc V_{1,m} + \sum_{m=1}^{2\bar {\mc E}}\mc V_{2,m},\\
	\mc V_{1,m}&=\int_{t-h_m}^t\tilde p(s)^\top S_m \tilde p(s)ds,\\
\mc V_{2,m}&= h_m\int_{t-h_m}^t (h_m+s-t)\dot{\tilde p}(s)^\top  R_m \dot{\tilde p}(s)ds.
\end{split}
\label{V2}
\endequ
where $V$ is  defined in \eqref{V} and $S_m\in\real^{(n-1)\times (n-1)}$ as well as $R_m\in\real^{(n-1)\times (n-1)}$ are positive definite matrices.
	
Recall that the proof of Proposition~\ref{stab} implies that with Assumption~\ref{ass:feas} there is a sufficiently small $\epsilon>0,$ such that $V$ is locally positive definite with respect to $z^*.$ Accordingly, for this value of $\epsilon>0$ also $\mc V$ is positive definite with respect to $z^*,$ since $\mc V_{1,m}$ is positive definite and so is $\mc V_{2,m}$, which can be seen in the following reformulation
\begin{equation*}
	\mc V_{2,m} = h_m\int_{-h_m}^0\int_{t+\phi}^t \left(\dot{\tilde p}(s)^\top  R_m \dot{\tilde p}(s)\right)dsd\phi.
\end{equation*}
Next, we inspect the time derivative of $\mc V$ along solutions of the system \eqref{eq: closed loop 6}. To this end, we at first assume $\epsilon=0.$ For that case by recalling \eqref{dotmcV1} together with the fact that 
\begin{equation}
\tilde p(t-\tau_m)=\tilde p(t)-\int_{t-\tau_m}^t \dot {\tilde p}(s)ds,
\label{eq:delay-fact}
\end{equation}
we have that
\begequ
\begin{split}
\dot V=&-\tilde \omega(t)^\top  D\tilde \omega(t) - \tilde p(t)^\top \sum_{m=1}^{2|\mc E_\ell|} \mc T_{\ell,m} \tilde p(t-\tau_m)\\
=& -\tilde \omega(t)^\top  D\tilde \omega(t) - \tilde p(t)^\top \sum_{m=1}^{2|\mc E_\ell|} \mc T_{\ell,m} \tilde p(t)\\
&+ \tilde p(t)^\top  \sum_{m=1}^{2|\mc E_\ell|} \mc T_{\ell,m}\int_{t-\tau_m}^t \dot {\tilde p}(s)ds\\
=&-\tilde \omega(t)^\top  D\tilde \omega(t) - \tilde p(t)^\top  \bar{\mc L}_{\ell}\tilde p(t) \\
&+ \tilde p(t)^\top  \sum_{m=1}^{2|\mc E_\ell|} \mc T_{\ell,m}\int_{t-\tau_m}^t \dot {\tilde p}(s)ds,
\end{split}
\label{dotmcV}
\endequ
where we have used \eqref{barLell} to obtain the last equality and recall that, for any $\ell\in\mc M,$ $\mc{\bar L}_\ell>0.$ Furthermore,
\begequ
	\dot {\mc V}_{1,m}= \tilde p(t)^\top  S_m  \tilde p(t)-\tilde p(t-h_m)^\top S_m  \tilde p(t-h_m).
	\label{dotV1}
\endequ
By using 
$$
\tilde p(t-h_m)=\tilde p(t)-\int_{t-h_m}^{t-\tau_m} \dot p(s)ds-\int_{t-\tau_m}^t \dot p(s)ds,
$$	
and defining
$$
\bar \eta_m=\col\left(\tilde p(t),\, \int_{t-h_{m}}^{t-\tau_{m}}\dot{\tilde p}(s) ds,\, \int_{t-\tau_{m}}^t \dot{\tilde p}(s)ds\right),
$$
\eqref{dotV1} is equivalent to	
\begequ
\begin{split}
\dot {\mc V}_{1,m}&=-\bar \eta_m^\top 
	\begin{bmatrix}
		\vectorzeros & -S_m & -S_m\\
		\ast  & S_m & S_m\\
		\ast  & \ast  & S_m	
	\end{bmatrix}
	\bar \eta_m.
\end{split}
\label{dotmcV_1}
\endequ		
Also, by differentiating $\mc V_{2,m}$ we obtain
\begequ
	\dot {\mc V}_{2,m}=-h_m\int_{t-h_m}^t\dot{\tilde p}(s)^\top  R_m \dot{\tilde p}(s)ds+h_m^2\dot{\tilde p}(t)^\top  R_m \dot{\tilde p}(t).
\label{dotmcV2}
\endequ
By following \cite{EF:14,EF:14b}, we reformulate the first term on the right-hand side of \eqref{dotmcV2} as
\begequ 
\begin{split}
&-h_m\int_{t-h_m}^t\dot{\tilde p}(s)^\top  R_m \dot{\tilde p}(s)ds=\\
&-h_m\int_{t-h_m}^{t-\tau_m}\dot{\tilde p}(s)^\top  R_m \dot{\tilde p}(s)ds-h_m\int_{t-\tau_m}^t\dot{\tilde p}(s)^\top R_m \dot{\tilde p}(s)ds
\end{split}
\label{jen1}
\endequ
Condition \eqref{rs12} is feasible by assumption. Thus, applying Jensen's inequality together with Lemma 1 in \cite{EF:14}, see also \cite{PP-JWK-CJ:11}, to both right-hand side terms in \eqref{jen1} yields the following estimate for the first term on the right-hand side of \eqref{dotmcV2}
\begequ
\begin{split}	
&-h_m\int_{t-h_m}^t\dot{\tilde p}(s)^\top R_m \dot{\tilde p}(s)ds\\
	&\leq-\begin{bmatrix}\int_{t-h_m}^{t-\tau_m}\dot{\tilde p}(s) ds\\\int_{t-\tau_m}^t \dot{\tilde p}(s)ds \end{bmatrix}^\top 
	\begin{bmatrix}
		R_m & S_{12,_m} \\\ast  & R_m
	\end{bmatrix}	
	\begin{bmatrix}\int_{t-h_m}^{t-\tau_m}\dot{\tilde p}(s) ds\\\int_{t-\tau_m}^t \dot{\tilde p}(s)ds \end{bmatrix}.
\end{split}
\label{dotmcV2_1}
\endequ
In order to rewrite the second term on the right-hand side of \eqref{dotmcV2}, i.e., $h_m^2\dot{\tilde p}(t)^\top  R_m \dot{\tilde p}(t),$ at first we replace $\dot{\tilde p}$ by its explicit vector field in \eqref{eq: closed loop 6}. This yields
\begequ
	\dot{\tilde p}(t)^\top  R_m \dot{\tilde p}(t)=\begin{bmatrix}
		\tilde \omega \\ \psi
	\end{bmatrix}^{\! \!\top} \!\!\!
	\begin{bmatrix} K^{\frac{1}{2}}WR_mW^\top K^{\frac{1}{2}} & -K^{\frac{1}{2}}W R_m\\
		-  R_m W^\top K^{\frac{1}{2}}&  R_m 
	\end{bmatrix}\!\!\!
	\begin{bmatrix}
		\tilde \omega \\\psi
	\end{bmatrix},
\label{dotmcV2_11}
\endequ
where $\psi = \sum_{j=1}^{2|\mc E|} \mc T_{\ell,j} \tilde p(t-\tau_j).$ 
Next, we apply the identity \eqref{eq:delay-fact} to obtain
\begequ
\psi=\sum_{j=1}^{2|\mc E|} \mc T_{\ell,j} \tilde p(t-\tau_j)=\sum_{j=1}^{2|\mc E_\ell|} \mc T_{\ell,j}\left(\tilde p(t)-\int_{t-\tau_j}^t \dot {\tilde p}(s)ds\right).
\notag
\endequ
Furthermore, we recall the fact \eqref{barLell} which implies that
$$
\sum_{j=1}^{2|\mc E_\ell|} \mc T_{\ell,j} \tilde p(t)=\mc{\bar L}_\ell\tilde p(t).
$$
Then, by introducing the short-hand vector
\begequ
\hat \eta=\col\left(
\tilde \omega,\, \tilde p(t),\, \sum_{j=1}^{2|\mc E_\ell|} \mc T_{\ell,j} \int_{t-\tau_j}^t \dot {\tilde p}(s)ds
\right),
\notag
\endequ
\eqref{dotmcV2_11} is equivalent to
\begequ
\begin{split}
	&\dot{\tilde p}(t)^\top  R_m \dot{\tilde p}(t)\\
	&=\hat \eta^\top
	\begin{bmatrix} K^{\frac{1}{2}}WR_mW^\top K^{\frac{1}{2}} & -K^{\frac{1}{2}}W R_m\bar{\mc L}_{\ell} & K^{\frac{1}{2}}W R_m \\
		\ast &  \bar{\mc L}_{\ell} R_m \bar{\mc L}_{\ell} & -\bar{\mc L}_{\ell} R_m \\
		\ast & \ast  & R_m
	\end{bmatrix}
	\hat \eta .
\end{split}
\label{dotmcV2_2}
\endequ
Finally, by collecting the terms \eqref{dotmcV}, \eqref{dotmcV_1}, \eqref{dotmcV2_1} and \eqref{dotmcV2_2}, $\dot {\mc V}$ can be upper-bounded by
\begequ
\dot{\mc V}\leq -\eta^\top  \Psi \eta,
\label{dotmcW}
\endequ
where $\eta \in\real^{(n+(1+4\bar{\mc E})(n-1))},$
\begequ
\begin{split}
\eta &= \col\left(\tilde \omega(t),\, \tilde p(t),\,\eta_1,\,\eta_2\right) ,\\
\eta_1&=\col\left(\int_{t-h_1}^{t-\tau_1}\dot{\tilde p}(s) ds,\, \ldots,\, \int_{t-h_{2\bar {\mc E}}}^{t-\tau_{2\bar {\mc E}}}\dot{\tilde p}(s) ds\right),\\
\eta_2&=\col\left(\int_{t-\tau_{1}}^t \dot{\tilde p}(s)ds ,\,\ldots ,\, \int_{t-\tau_{2\bar {\mc E}}}^t \dot{\tilde p}(s)ds\right)
\end{split}
\endequ
and $\Psi$ is defined in \eqref{psi}. 
Therefore, if conditions \eqref{psi}, \eqref{rs12} are satisfied, then $\dot {\mc V}\leq0.$

 As of now, $\epsilon=0$ and $\dot{\mc V}$ is not strict, i.e., not negative definite in all state variables. Yet, as the system \eqref{eq: closed loop 6} is non-autonomous, we need a strict common LKF in order to establish the asymptotic stability claim. By making use of Proposition~\ref{stab}, this can be achieved without major difficulties as follows.
		From \eqref{dotmcV1} in the proof of Proposition~\ref{stab} together with \eqref{dotmcW}, we have that for $\epsilon\neq0,$ 
		\begequ
		\begin{split}
			\dot {\mc V}&\leq-\bar\xi^\top \left(  \begin{bmatrix} \vectorzeros[n\times n] & \vectorzeros[n \times (n+(1+4\bar{\mc E})(n-1))]\\
			\ast  & \Psi\end{bmatrix} + \epsilon  \Xi \right)
			\bar\xi,
		\end{split}
		\label{dotmcW2}
		\endequ
		where
		\begequ
		\bar\xi\!=\!\col\!\left(\!\nabla_{\tilde \theta} U(\tilde \theta(t)+\theta^*)\!-\!\nabla_{\tilde \theta} U(\theta^*),\tilde \omega(t), \, \tilde p(t), \eta_1, \, \eta_2\right)\!,
		\notag
		\endequ
		\begequ
		\Xi=\begin{bmatrix} A & 0.5 AD & 0.5  AK^{\frac{1}{2}}W & \vectorzeros& \vectorzeros\\
		\ast  & -0.5 E_{22} & \vectorzeros &\vectorzeros& \vectorzeros\\
		\ast  & \ast & \vectorzeros & \vectorzeros& \vectorzeros\\
		\ast & \ast & \ast  & \vectorzeros& \vectorzeros\\
		\ast & \ast & \ast  & \ast& \vectorzeros\\
		\end{bmatrix},
		\notag
		\endequ
		$E_{22}$ is defined in \eqref{e22}
		and $\vectorzeros$ denotes a zero matrix of appropriate dimensions.
		Recall that if conditions \eqref{psi}, \eqref{rs12} are satisfied, then $\Psi >0$ for all $\bar{\mc L}_{\ell}$, $\ell\in\mc M.$
		Thus, following the proof of Proposition~\ref{stab}, Lemma \ref{Lemma: matrix regularization} in Appendix~\ref{appendix2} implies that there is a sufficiently small $\epsilon>0$ such that the matrix sum on the right-hand side of \eqref{dotmcW2} is positive definite and, at the same time, $\mc V$ is locally positive definite. Consequently, there exists a real $\gamma>0$ such that
		$
		\dot{\mc V}(t)\leq -\gamma\|z(t)\|,$ $z(t)=\col(\tilde \theta(t),\,\tilde\omega(t),\,\tilde p(t)).$
		Local uniform asymptotic stability of $z^*$ follows by invoking the standard Lyapunov-Krasovskii theorem \cite{EF:14,EF:14b} and arguments from \cite{EF-AS-JPR:04} for systems with piecewise-conti\-nuous delays. 
		By direct inspection, if \eqref{psi}, \eqref{rs12} are satisfied for some $h_m \in \real_{\geq0},$ then they are also satisfied for any $\tau_m(t) \in [0,h_m),$ which in particular includes the case of constant delays, i.e., $\dot \tau_m=0.$ This completes the proof.

\end{document}